\documentclass[12pt]{article}
\usepackage{amssymb, amscd, amsmath}

\newtheorem{theorem}{Theorem}
\newtheorem{lemma}[theorem]{Lemma}

\newtheorem{proposition}[theorem]{Proposition}

\newcommand{\C}{\kappa}

\newcommand{\jj}{\mathcal{J}} 
\newcommand{\E}{\mathcal{E}}

\newcommand{\sym}{\mathrm{Sym}} 

\newcommand{\con}{N^\vee} 
\newcommand{\spec}{\mathrm{Spec}} 
\newcommand{\hato}{\widehat{\mathcal{O}}} 
\newcommand{\ox}{\mathcal{O}_X} 

\begin{document}

\title{Obstructions to extension of vector bundles}
\date{September 30, 2022}

\author{Vladimir Baranovsky, Hongseok Chang}

\maketitle

\abstract{In the holomorphic or algebraic setting we consider 
a vector bundle $E$ on a smooth subvariety $X$ in a smooth
variety $Y$ over a field of characteristic zero. Assuming 
$E$ extends to the $l$-th formal neighborhood of $X$ in $Y$,
we study cohomological obstructions to extending it
further to the $k$-th neighborhood,  for $k > l$.}

\section*{Introduction}

We consider the setup of a smooth closed 
algebraic subvariety $X$ of dimension $p$, in a smooth 
algebraic variety $Y$ of dimension $p+q$ over a field of characteristic zero 
$\kappa$, and a vector bundle $E$ of rank $e$ on $X$. 
The $k$-th neighborhood $Y^{(k)}$ of $X$ in $Y$ is the scheme with the structure sheaf 
$\mathcal{O}^{(k)} = \mathcal{O}_Y/I^{k+1}_X$ where $\mathcal{O}_Y$ is the structure sheaf
of $Y$ and $I_X$ is the ideal sheaf of $X$ (and the underlying topological space of 
$Y^{(k)}$ is homeomorphic to that of $X$). 

We are interested in the question of 
existence of a locally free sheaf $E^{(k)}$ of $\mathcal{O}^{(k)}$-modules, 
such that $E^{(k)}/I_X E^{(k)}$ is isomorphic to $E$ as a sheaf of 
$\mathcal{O}_X$-modules. We will call such an $E^{(k)}$ an \textit{order $k$
extension} of $E$.
In a similar $C^\infty$ situation the problem would be 
trivial due to existence of a tubular neighborhood $U \subset Y$ of $X$ which 
admits a projection $U \to X$ so we can extend $E$ even to $U$ by taking the 
pullback with respect to the projection. 

In the algebraic situation though we only have a formal analogue of the 
tubular neighborhood theorem valied when $Y$ and $X$ are affine: in this case
the limit $\mathcal{O}^{(\infty)} = \lim_k \mathcal{O}^{(k)}$ is
non-canonically isomorphic to the completed symmetric algebra of the conormal 
bundle $N^*$ of $X$ in $Y$. As the isomorphism is non-canonical, it will not 
hold for general $(X, Y)$ and instead we have a twisted version of the statement, 
in which $\mathcal{O}^{(\infty)}$ is described via a structure of an 
$L_\infty$-algebroid on the shift $N[-1]$ of the normal bundle. 

This makes our question about existence of $E^{(k)}$ non-trivial. 
We apply the standard machinery of formal geometry and obstruction theory 
to find cohomology classes that must vanish for a choice of $E^{(l)}$ to 
admit an extension to some $E^{(k)}$. For $k \leq 2l+1$ this is also 
a sufficient condition while for $k \geq 2l+2$ one is dealing with 
abelian "shadows" of non-abelian cohomology classes. 

In more detail, for $k \leq 2l$ we find a cohomology class in 
$H^2(Y^{(l)}, (\jj^{l+1}/\jj^{k+1}) \otimes_{\mathcal{O}^{(l)}} E^{(l)})$
(where $\jj \subset \mathcal{O}^{(\infty)}$ is the ideal sheaf of $X$, 
i.e. completion of $I_X$) 
that must vanish if $E^{(k)}$ exists. If the vanishing does hold
then the set of isomorphism classes of such bundle extensions is
a nonempty torsor over 
$H^1(Y^{(l)}, (\jj^{l+1}/\jj^{k+1}) \otimes_{\mathcal{O}^{(l)}} E^{(l)})$.

We describe three different approaches to this problem. The first is
based on formal geometry, works over a field of characterstic zero. 
It involved infinite dimensional torsors which makes it perhaps less
explicit than the other two approaches. The second approach works
over $\kappa = \mathbb{C}$ and involves the usual machinery of connections 
and similar operators. The third approach involves affine open covers
and Cech complexes, although we only make it explicit for $k=1, 2$
(this is the situation when a relevant Cech complex has a structure 
of a dg Lie algebra, rather than an $L_\infty$ algebra). 

\bigskip
\noindent
\textbf{Remarks.}

\begin{enumerate}
\item Although we will not use this, but constructing an infinite
order extension $E^{(\infty)}$ is equivalent to giving $E$ the structure
of a module over the $L_\infty$-algebroid, i.e. an 
$L_\infty$-morphism $N[-1] \to At(E)$ into the Atiyah algebra sheaf
of $E$, such that the composition with the canonical Lie morphism 
$At(E) \to T_X$ equals the anchor $L_\infty$-morphism of $N[-1]$. 
In fact, the operators $\mathfrak{M}_i$ in Section 4 below are adjoint
to the components of such an $L_\infty$ morphism. 

\item A special case of this situation, when $Y = X \times X$
and $X \to Y$ is the diagonal morphism, was studied in \cite{Ka}. In this
case $E$ admits an extension to $Y$ by pullback with respect to either
of the two projections, hence and bundle $E$ admits a structure
of an $L^\infty$-module over $N[-1] \simeq T[-1]$, as discussed in 
Section 2.7 of \cite{Ka}. 
Note that the anchor morphism is trivial in this case. 
\end{enumerate}

\bigskip
\noindent 
This paper is organized as follows. In Section 1 we consider the
affine case and show that the pair $(\mathcal{O}^{(\infty)}, 
E^{(k)})$ is can be trivialized for any choice of $E^{(k)}$.
Section 2 outlines our strategies for dealing with the general 
case, as long as a simple computation involving lift of
Maurer-Cartan elements, on which most identities of the paper
are based. In Section 3 we outline the basic constructions of 
formal geometry, applications to lifts of torsors and do Lie
cocycle computations to make our obstruction theory statements 
more explicit. In Section 4 we discuss the Dolbeault approach 
and, in particular, Yu's formulation for the data describing 
$\mathcal{O}^{(\infty)}$, and give our description for the 
data describing $E^{(k)}$. Finally, in Section 5, we give
a Cech version of $k=1,2$ and briefly indicate how the approach 
can be extended to arbitrary $k > 0$.

\bigskip
\noindent
\textbf{Notation.} We will suppress $\mathcal{O}_X$ from notation in 
expressions like $\otimes_{\mathcal{O}_X}$, $Hom_{\mathcal{O}_X}$, 
$End_{\mathcal{O}_X}$
or $Sym^s_{\mathcal{O}_X}$, while keeping the reference to the base 
ring in all other cases, such as $\otimes_{\kappa}$.

\section{Formal completion: the affine case.}

Here we discuss the case of a closed embedding
of smooth affine varieties

\begin{lemma} 
\label{lemma: filtered k isomorphism}
    Given two smooth affine varieties
    $        X = \spec \, R/I \hookrightarrow Y = \spec \, R$
    there exists a filtered $\C$-algebra isomorphism between
    \begin{equation} \label{filtered k isomorphism}
        \lim_{\longleftarrow} R/I^{u+1} \quad \simeq \quad \widehat{\sym}_R^\bullet \con = \prod_{u \geq 0} \sym_R^u \con
    \end{equation}
    where $\con:= I/I^2$ is the conormal module. This isomorphism may be assumed to 
    induce the canonical isomorphism on associated graded factors: $I^u/I^{u+1}   \simeq \sym^u_R \con$. 
    
    \medskip
    \noindent 
    Moreover, given a projective $R/I$-module $E$
    and a projective $R/I^{k+1}$-module $E^{(k)}$ with an isomorphism of 
    $R/I$-modules     $E  \simeq E^{(k)}/I E^{(k)}$ , there exists an isomorphism 
    \begin{equation} \label{lemma:filtered bundle isomorphism eq:filtered bundle isomorphism}
    	E^{(k)} \simeq \big(\prod_{u = 0}^k \sym_R^u \con\big) \otimes_R E
    \end{equation}
 of filtered $\kappa$-vector spaces, which is compatible with the module actions of the filtered $K$
 -algebras in (\ref{filtered k isomorphism})
\end{lemma}

\medskip
\noindent
\textit{Proof:} The isomorphism of filtered algebras is proved, e.g., in Lemma 5.2 of \cite{CCT}. 
Once the filtered algebras are identified we can view 
$\sym^{u}_R \con$ as a subset of $R/I^{k + 1}$. 
For the case of projective modules we proceed by induction on $k \geq 0$.
For $k = 0$ there is nothing to prove.  Assuming the isomorphism exists for a particular $k > 0$, observe that
$$
I^{k+1} E^{(k+1)} /I^{k + 2} E^{(k+1)} \simeq  
I^{k + 1} /I^{k + 2}  \otimes_{R/I^{k + 2} } E^{(k+1)}
	\simeq \sym_R^{k +1} \con \otimes_R E 
$$
Hence we have an exact sequence of $R/I^{k+1}$ modules
$$
0 \to \sym_R^{k +1} \con \otimes_R E  \to  E^{(k+1)} \to 
 E^{(k)} \to 0
$$
and by inductive assumption the quotient is isomorphic to 
$\big(\prod_{u = 0}^k \sym_R^u \con\big) \otimes_R E$. Thus we can 
choose a vector space splitting 
$$
 \Phi: \big(\prod_{u = 0}^{k+1} \sym_R^u \con\big) \otimes_R E 
\simeq  E^{(k+1)}
$$
and set for  
$x \in R^{\leq k+1}  = \prod_{u = 0}^{k+1} \sym_R^u \con$
and $m \in  R^{\leq k} \otimes_R E $
$$
\psi(x, m) = x \Phi(m) - \Phi (xm) 
\in \sym_R^{k +1} \con \otimes_R E. 
$$
It is immediate that the Hochschild cochain condition 
$\psi(xy, m) - x \psi(y, m) - \psi(x, ym) =0$ is satisfied. By Lemma 9.1.9
in \cite{We} the 
Hochschild complex of $R^{\leq k+1}$ with values in the bimodule
$Hom_\kappa( R^{\leq k} \otimes_R E ,
\sym_R^{k +1} \con \otimes_R E)$ computes 
$Ext^\bullet_{R^{\leq k + 1}} (R^{\leq k} \otimes_R E, \sym_R^{k +1} \con \otimes_R E)$. 
Using the short exact sequence 
$$
0 \to \sym_R^{k +1} \con \otimes_R E \to 
R^{\leq k+1} \otimes_R E \to R^{\leq k} \otimes_R E \to 0
$$
and the fact that the natural pullback maps give isomorphisms
$$
Ext^0_{R^{\leq k + 1}} (R^{\leq k} \otimes_R E, \sym_R^{\leq k + 1} \con \otimes_R E) \simeq Ext^0_{R^{\leq k + 1}} (E, \sym_R^{k +1} \con \otimes_R E) \simeq 
$$
$$
\simeq 
Ext^0_{R^{\leq k + 1}} (R^{\leq k + 1} \otimes_R E, \sym_R^{k +1} \con \otimes_R E)
$$
we conclude that 
\begin{equation}
\label{ext-hom}
Ext^1_{R^{\leq k + 1}} (R^{\leq k} \otimes_R E, \sym_R^{k +1} \con \otimes_R E)
\simeq Hom_{R^{\leq k + 1}} (\sym_R^{\leq k + 1} \con \otimes_R E, \sym_R^{\leq k + 1} \con \otimes_R E)
\end{equation}
Explicitly, the isomorphism is induced by restricting a Hochschild cocycle to 
$\sym_R^{k +1} \con \subset R^{\leq k + 1}$ and then observing that this restriction 
factors through the surjection to $E$. In particular,
$E^{(k +1)}$ is isomorphic to 
 $\big(\prod_{u = 0}^{\kappa + 1} \sym_R^u \con\big) \otimes_R E $ via a filtered isomorphism if both define the same class in the 
 left hand side of \eqref{ext-hom}.  But both 
$\psi$ and a similar element for 
 $\big(\prod_{u = 0}^{k + 1} \sym_R^u \con\big) \otimes_R E $ give the identity 
 element on the right hand side of \eqref{ext-hom}, finishing the proof.  $\square$

\section{General strategy}

\subsection{Three approaches}
In the general case the isomorphism proved in Lemma 
\ref{lemma: filtered k isomorphism} does not exist. 
This can be remediated in three different ways

\begin{enumerate}
	
	\item \textit{Formal geometry.} We can replace $X$ by a torsor $P(k) \to X$ over 
	a prounipotent group, understood as a generalized cover of $X$.  The fiber of $P(k)$
	over a point $x \in X$ parameterizes isomorphisms of formal completions at the point $x$ 
	$$
	\widehat{\mathcal{O}}_{X, x} \simeq \mathcal{A}, \quad
	\widehat{\mathcal{O}}^{\infty}_x  \simeq \prod_{u = 0}^{\infty}
	Sym^u_{\mathcal{A}} \mathcal{S}, \quad 
	\widehat{E}^{(k)}_x \simeq \bigoplus_{u=0}^k  Sym^u_{\mathcal{A}}
	\mathcal{S} \otimes_{\mathcal{A}} \E 
	$$ 
	with the ``standard formal models" built out of a power series algebra $\mathcal{A}$ and 
	two projective (hence free) modules $\mathcal{S}$, $\mathcal{A}$ over it. $P(k)$ has a structure of a torsor over infinite dimensional prounipotent algebrac group $G(k)$. In the torsor is non-trivial, 
	so a global algebraic section does not exist.  Instead
	of working on $X$ we work with $G(k)$-equivariant objects on $P(k)$, or a certain quotient space $Q(k)$ of it which is torsor over a sheaf
	of unipotent groups on $X$. 
	
		\item \textit{Dolbeault model}. When 
		$\C = \mathbb{C}$ we can choose a $C^\infty$-section 	
		$s: X \to Q(k)$ and consider sections various bundles 
		and sections that reflect deviation of $s$ from being  holomorphic. 

	\item \textit{Cech model}. 
	We can choose an open cover $X = \bigcup U_i$ and then local sections 
	$s_i: U_i \to Q(k)|_{U_i}$, paying attention to what happens on double and triple
	intersections.

\end{enumerate}

\subsection{Reminder on lift of Maurer-Cartan elements.}

Eventually, all three approaches to extending the bundle $E$ to from an order $k$  to 
an order $k+1$  will rely on 
a fairly easy formalism of lifting Maurer-Cartan solutions, which 
we explain below to fix notation. 
Consider an abelian extension of differential graded Lie algebras
$$
0 \to \mathcal{E} \to \widehat{\mathcal{L}} \to \mathcal{L} \to 0.
$$
where $\mathcal{E}$ is a dg module over $\mathcal{L}$. 
We are interested in the question of lifting a (degree one) 
Maurer-Cartan element in 
$\varphi \in \mathcal{L}^1$ to  $\widehat{\mathcal{L}}^1$.  
Suppose we have a linear section $s: \mathcal{L} \to  \widehat{\mathcal{L}}$ which we
do not assume to commute with differential or bracket. Thus  we have 
possibly nontrivial maps of degree 1 and 0, respectively
$$
\Delta_1: \mathcal{L} \to  \mathcal{E};
\qquad 
\Delta_1(x) = (d_{\widehat{\mathcal{L}}} s - s d_{\mathcal{L}})(x),
$$
$$
\Delta_2: \Lambda^2\mathcal{L} \to  \mathcal{E} ;  \qquad
\Delta_2(x_1 \wedge x_2) = [s(x_1), s(x_2)]_{\widehat{\mathcal{L}}} - s[x_1, x_2]_{\mathcal{L}}
$$
Now suppose that $\varphi$ is lifted to a Mauer-Cartan solution
$
\widehat{\varphi} = s(\varphi) + \alpha \in \mathcal{L}^1 
$
with $\alpha \in \mathcal{E}^1$.  This leads to the equations 
$$
0 = d (s\varphi + \alpha)  + \frac{1}{2}[s\varphi + \alpha, s\varphi + \alpha]  =
$$
$$
= s(d(\varphi) + \frac{1}{2} [\varphi, \varphi]) + \Delta_1(\varphi) + d\alpha + 
 [s\varphi, \alpha]
+ \frac{1}{2} \Delta_2(\varphi, \varphi)
$$
In other words, using the Maurer-Cartan equation for $\varphi$ 
and the fact that the commutator with $s \varphi$ agrees with 
the action of $\mathcal{L}$ on $\mathcal{E}$ we can write 
\begin{equation}
	\label{lift_equations}
	-(d+ \varphi \cdot)\alpha = \Delta_1(\varphi)  + \frac{1}{2} \Delta_2(\varphi, \varphi)  
\end{equation} 
This can be reformulated in the usual way: a choice of $\alpha$ is possible if  the right hand side of the 
second equation in (\ref{lift_equations}) gives a zero 
class in the cohomology of $(\mathcal{E}, d + 
\varphi \cdot) $.

\section{The formal geometry model}

\subsection{Torsors and Lie cohomology classes.}

We give here a short introduction to  Formal Geometry, a specific version needed for our purposes.
Suppose that $G$ is a proalgebraic group which splits into a semidirect product 
$G \simeq U \rtimes F$ of a prouniponent normal subgroup $U$ and a group 
$F$.  Let $P \to X$ be a torsor over $G$ (in our case, in the Zariski topology) 
and $\pi: Q = P/G \to X$. 
For any coherent $\mathcal{O}_X$-module 
$M$ its pullback $\pi^* M$ has a canonical flat connection along the fibers 
and the relative de Rham  complex $\pi^* M \otimes_{\mathcal{O}_X} 
\Omega^\bullet_{\pi}$ is well defined. Since $Q$ is a limit of 
fibrations with finite dimensional
affine fibers, the canonical map 
$$
M \to \pi_* \big(\pi^* M \otimes_{\mathcal{O}_Q} 
\Omega^\bullet_{\pi} \big) 
$$
is a quasi-isomoprhism. 
A slighly more general setting is to enlarge $\mathfrak{g}$ to a 
Lie algebra $\mathfrak{h}$ in such a way that $(G, \mathfrak{h})$
forms a Harish-Chandra pair, see Section 2 of \cite{BK}. 
If $P$ has an additional transitive Harish-Chandra structure, cf.
\cite{BK}, and the pullback of $M$ to $P$ is isomorphic to 
a trivial vector bundle $\mathcal{O}_P \otimes \kappa \mathcal{M}$
associated with a module $\mathcal{M}$ over $(G, \mathfrak{h})$
then $M$ has a flat connection and its de Rham complex is resolved
by the pushforward of the de Rham complex of 
$\mathcal{O}_P \otimes \kappa \mathcal{M}$ to $X$. 
The main source of examples is when $M$ is a jet bundle of a certain 
vector bundle $L$ and its local de Rham cohomology is only nontrivial 
in degree 0 where it is equal to $L$. In this case the global 
de Rham cohomology of $M$ is isomorphic to coherent cohomology of $L$.

We will consider the complex of relative Lie cochains on the algebra $\mathfrak{g} = Lie(G)$, 
with respect to the Lie subalgebra $\mathfrak{f} = Lie(F)$, with coefficients in a
$\mathfrak{g}$-module $\mathcal{M}$.  
In degree $j$ its component
$C^j(\mathfrak{g}, \mathfrak{f}; \mathcal{M})$ is formed by 
$\mathfrak{f}$-invariant
cochains $\alpha: \Lambda^j (\mathfrak{g}/\mathfrak{f}) \to \mathcal{M}$ with the Lie cochain differential $\delta_{Lie}$ as in Section 1.3 of 
\cite{Fu}. 

Since $\mathfrak{g}$ can be identified with the fiber of relative tangent bundle for 
$P \to X$, a Lie cochain gives a differential form on $P$ with values in $\mathcal{O}_P \otimes_\kappa
\mathcal{M}$, and relative cochains may be identified with forms on $Q$ with values in 
$\mathcal{O}_Q \otimes_{\kappa} \mathcal{M}$. The following construction is described
in slightly different terms in \cite{Fu}. 

\begin{lemma}
	Let $M$ be a locally free sheaf on $X$ and $\mathcal{M}$ a module over $G$ such that
	the pullback of $M$ to $P$ is isomorphic to $\mathcal{O}_P \otimes_{\C} \mathcal{M}$
	as a $G$-equivariant $\mathcal{O}_P$-module. 
	
	Let $\alpha \in C^j(\mathfrak{g}, \mathfrak{f}; \mathcal{M})$ be a relative Lie cochain.
	Then $\alpha$ defines a global section of $ \pi^* M \otimes_{\mathcal{O}_P} 
	\Omega^j_{\pi}$ and hence of $\pi_* \big(\pi^* M \otimes_{\mathcal{O}_P} 
	\Omega^\bullet_{\pi} \big) $. If $\alpha$ is closed it defines 
	the \emph{Gelfand-Fuks cohomology class}
	$GF(P, \alpha) \in H^j(X, M)$ (depending on the choice of $P$).
	
	In the Harish-Chandra setting, for $\alpha \in C^j(\mathfrak{h}, \mathfrak{f}, 
	\mathcal{M})$ the class $GF(P, \alpha)$ is defined 
	in the de Rham cohomology $H^j_{DR}(X, M)$.
\end{lemma}

\bigskip
\noindent 
Now suppose we have a morphism of pro-algebraic groups $\widetilde{G} \to G$ 
with a unipotent splitting as above, and assume that it induces an isomorphism 
$\widetilde{F} \to F$ and a morphism of unipotent parts $\widetilde{U} \to U$.
If a $G$-torsor $P \to X$ lifts to a $\widetilde{G}$-torsor $\widetilde{P} \to P 
\to X$, then the induced projection $\widetilde{\pi}: \widetilde{Q} = \widetilde{P}/F \to 
X$ factors through $\pi$ which allows to compare the two resolutions of $M$ on $X$:
$$
M \to \pi_* (\pi^* M \otimes_{\mathcal{O}_Q} \Omega^\bullet_{\pi}) 
\to \widetilde{\pi}_* (\widetilde{\pi}^* M \otimes_{\mathcal{O}_{\widetilde{Q}}} 
\Omega^\bullet_{\widetilde{\pi}}) 
$$

\begin{lemma}
\label{torsor extension}
	Let $\alpha \in C^j(\mathfrak{g}, \mathfrak{f}; \mathcal{M})$ be
	a relative Lie cocycle that has an exact pullback  
	$\widetilde{\alpha} \in C^j(\widetilde{\mathfrak{g}}, \mathfrak{f}; \mathcal{M})$. 
	Then the cohomology class of $GF(P, \alpha)$ is zero.  If $\widetilde{G}$ is an abelian extension 
	that splits over $F \subset G$
	$$
	1 \to \mathcal{M} \to \widetilde{P} \to P \to 1
	$$
	with $\mathcal{M}$ viewed as abelian group; 
	and $\alpha \in C^2(\mathfrak{g}, f; \mathcal{M})$ is the extension 
	cocycle for the induced extension of Lie algebras, 
	 then the vanishing of $GF(P, \alpha)$ is equivalent to existence of 
	$\widetilde{P}$.  	Similar statements hold in Harish-Chandra setting. 
\end{lemma}
\textit{Proof.} The first statement follows from the 
comparison of resolutions before the lemma. The second statement is
proved in Proposition 2.7 of \cite{BK} $\square$

\subsection{Formal derivations and cocycles.}

We now proceed to discuss the Lie algebras to which the above theory will be applied, 
delaying the description of the proalgebraic groups and their torsors until later. 
Let
$\mathcal{A} = \mathbb{C}[[x_1, \ldots, x_p]]$ and let $\mathcal{S}$, 
$\mathcal{E}$ be fixed free $\mathcal{A}$-modules of ranks $q, e$, respectively. 
Define
$$
\mathcal{S}^u = Sym^u_{\mathcal{A}} \mathcal{S}, \quad 
\mathcal{A}_S = \prod_{u \geq 0} \mathcal{S}^u, \quad \E^{\leq k} = \bigoplus_{u=0}^k 
\mathcal{S}^u \otimes_{\mathcal{A}} \E,
$$
viewing $\E^{\leq k}$ as a module over $\mathcal{A}_\mathcal{S}$.
Consider the Lie algebra of filtration-preserving pair derivations of 
$$
Der_k = Der(\mathcal{A}_\mathcal{S}, \E^{\leq k}), 
$$
i.e. pairs of maps $\varphi: \mathcal{A}_\mathcal{S} \to \mathcal{A}_\mathcal{S}$, 
$\psi: \E^{\leq k} \to \E^{\leq k}$ such that $\varphi$ is a derivation and 
$$
\psi(a \cdot e) = \varphi(a) \cdot e + a \cdot \psi(e).
$$
We require that $\varphi$ is continuous in the topology 
defined by the powers of maximal 
ideal $\mathfrak{m} \subset \mathcal{A}$ generated by the variables 
$x_1, \ldots, x_p$, and also the filtration by powers of the ideal in $\mathcal{A}_\mathcal{S}$ 
generated by $\mathcal{S}$. Similar condition is imposed on $\psi$. 
Then $Der_{-1} =: Der(\mathcal{A})_{\mathcal{S}}$ is the Lie algebra of derivations $\varphi$ satisfying the 
above properties.

Due to the Leibniz rule, such a pair can be described by operators 
$A_v: \mathcal{A} \to \mathcal{S}^v$, 
$L_v: \mathcal{S}=\mathcal{S}^1 \to \mathcal{S}^{v+1}$, $M_v: \E \to \mathcal{S}^v \otimes_{\mathcal{A}} \E$ for $v \geq 0$.
Here each $A_v$ is a derivation, which can also be described by an 
$\mathcal{A}$-linear map $a_v: \Omega^1_{\mathcal{A}} \to \mathcal{S}^v$, and 
$L_v$, $M_v$  are fist order operators with scalar symbols given by 
$a_v$ (composed with the symmetric algebra product). 
We can also write 
$$
Der_k = \prod_{u \geq 0} Der^u_k
$$
where $Der^u_k$ is spanned by 
triples $(A_u, L_u, M_u)$ with fixed $u$. Defining 
$Der^+_k$ as the infinite product $\prod_{u \geq 1} Der^u_k$, have a
semidirect product splitting. 
$$
Der_k \simeq Der^+_k \rtimes Der^0_k, 
$$

\bigskip
\noindent 
We would like to focus on the extension of Lie algebras
\begin{equation}
	\label{main_extension}
0 \to End^{l, k} = \bigoplus_{u = l+1}^k  End^u \to Der_k \to Der_l \to 0, \qquad End^u: = 
\mathcal{S}^u \otimes_{\mathcal{A}} End_{\mathcal{A}}(\E)
\end{equation}
for $k > l$.  This extension does split into a semi-direct product, but we would like to 
make the splitting more explicit. 
Choose and fix a connection $\nabla^{\E}$ on $\E$ which exists since $E$ is free
over $\mathcal{A}$. A derivation $\varphi_{v}: \bigoplus_{u \geq 0} \mathcal{S}^u \to \bigoplus_{u \geq 0} \mathcal{S}^{u+v}$
defines operators 
$
a_v \nabla^\E: \E \to \mathcal{S}^v \otimes_\mathcal{A} \E
$
which we extend to 
$$
 \qquad s\varphi_v: \bigoplus_{u = 0}^{k-v} \mathcal{S}^u 
\otimes_{\mathcal{A}} \E \to \bigoplus_{u=0}^{k-v} \mathcal{S}^{u+v} \otimes_{\mathcal{A}} \E,
 \quad u = 0, \ldots, k 
$$
Leibniz rule and 
restriction of $\varphi_v: \mathcal{S}^u \to \mathcal{S}^{u+v}$. 
Summing over $v = 0, \ldots, k$
we obtain
a derivation of $\E^{\leq k}$ that we denote by $s\varphi$.
A general element in $Der_k$ will have components $\psi_v = s\varphi_v
+ e_v$ where $e_0, \ldots, e_k$ are $\mathcal{A}$-linear.  
When $k > l$ we use the same notation $s$ for a linear splitting 
$s: Der_l \to Der_k$ which takes 
$$
\big((\varphi_0, \varphi_1, \ldots), (s\varphi_0 + e_0, \ldots, s\varphi_l + e_l)\big) \mapsto
\big((\varphi_0, \varphi_1, \ldots), (s\varphi_0 + e_0, \ldots, s\varphi_l + e_l, 
s\varphi_{l+1}, \ldots, s\varphi_k)\big)
$$
\begin{lemma}
	The splitting $s: Der_l \to  Der_k$  induced by $\nabla^{\E}$ is 
	compatible with the Lie bracket if and only if $\nabla^{\E}$ is flat.
\end{lemma}
\textit{Proof.} First assume that the splitting is compatible with the 
bracket. Restricting to the degree 0 part, we get a Lie splitting of 
the Atiyah Lie algebra $Der(\mathcal{A}, \E)$ of filtration preserving 
derivations of the pair $(\mathcal{A}, \E)$
$$
0 \to End^0 \to Der(\mathcal{A}, \E) \to T_{\mathcal{A}} \to 0
$$
given by the connection $\nabla^{\E}$ and its compatibility with the 
bracket is equivalent to the vanishing of the curvature by 
definition of curvature. 

For the other direction, it suffices to assume $l=-1$.  
Let $\E_0$ be the quotient of $\E$ modulo the maximal ideal of
$\mathcal{A}$.
As $\nabla^\E$ is flat, there is 
a unique section of $\E \to \E_0$ with image in the subspace of flat sections. 
Moreover, this gives an
isomorphism of $\mathcal{A}$-modules
\begin{equation}
	\label{flat_splitting} 
\E \simeq \E_0 \otimes_{\kappa} \mathcal{A} 
\end{equation}
such that $\nabla^\E$ corresponds to the usual vector field action on the second tensor
factor. 
Then this isomorphism sends $\nabla^\E(\varphi) $to $Id_{\E_0} \otimes \varphi$ 
and the lemma follows. 
$\square$

\bigskip
\noindent 
The Lie algebra extension \eqref{main_extension} will have abelian kernel if
$k \leq 2l +1$. In full generality, we can  obtain a quotient extension that has 
abelian kernel. To give more detail, write $End_0^j \subset End^j$ for the subalgebra
of endomorphisms of trace zero (twisted by appropriate symmetric power of $\mathcal{S}$). 
Then, for $k > 2l +1$, the subalgebra spanned by the commutators in $End^{l, k}$ is 
$$
[End^{l, k}, End^{l, k}] = \bigoplus_{i = 2l+2}^k End_0^i 
$$
It is easy to see that this is an ideal in $Der_k$ and, denoting 
$$
End_{ab}^{l,k} : = End^{l, k}/ [End^{l, k}, End^{l, k}] 
\simeq \bigoplus_{i = l+1}^{2l+1} End^i  \oplus \mathcal{S}^{2l+2} \oplus \ldots 
\oplus \mathcal{S}^k  
$$
we obtain an abelian extension (that splits into a semidirect product, once a flat connection on 
$\E$ is chosen):
\begin{equation}
	\label{extension_ab}
0 \to End_{ab}^{l, k} \to Der_k' \to Der_l \to 0
\end{equation}

\bigskip
\noindent 
Fix a flat connection $\nabla^\E$ on $\E$. 
The Lie cocycle $c(l,k)$ corresponding to the abelian extension \eqref{extension_ab}
has homogeneous components in degree $i \in \{l+1, \ldots, k\}$.   It
measures incompatibily of the section with the Lie bracket in $Der'_k$.
For $l+1 \leq v \leq 2l+1$ the component is 
\begin{equation}
\label{cocycle}
c(l,k)_v ((\varphi, \psi) \wedge (\widetilde{\varphi}, \widetilde{\psi})) = 
\sum_{p = 0}^{l} \Big([s\varphi_{v-p}, \widetilde{e}_p] + 
[s\widetilde{\varphi}_{v-p}, e_p]\Big)+ 
\sum_{p = v-l}^{l} [e_{v-p}, \widetilde{e}_p]
\end{equation}
We observe that the second sum
has $l$ terms for $v = l+1$ and 0 terms for $v = 2l+1$. 
For $v \in \{2l+2, \ldots, k\}$ (if that range is non-empty) the second sum is trivial 
and the first is replaced by its trace. 
The simplest statement can be made under the flatness assumption.  

\begin{proposition} 
\label{proposition: cocycle}
	If the splitting of 
	\eqref{extension_ab} is induced by a flat connection $\nabla^\E$, the extension cocycle $c(k,l)$
	of \eqref{extension_ab} is relative with respect to the subalgebra $Der'_0 \subset Der_0$ formed by all
	derivations $(\varphi_0+\varphi_1+\ldots, s\varphi_0 + e_0)$ with constant $e_0$ (which is well defined by flatness, cf. \eqref{flat_splitting}). 
	Its image  in $C^2(Der'_k, Der_0'; End_{ab}^{l, k}) $ is equal to $-\delta \beta$
	where $\beta: Der'_k \to End_{ab}^{l, k}$ is the projection (viewed as 1-cochain) induced by 
	the splitting of 	\eqref{extension_ab}.	
	In particular, the cohomology class
	$[\alpha]$ is in the kernel of the map 
	$$
	H^2(Der_l, Der_0'; End_{ab}^{l, k}) \to H^2(Der'_k, Der_0'; End_{ab}^{l, k}) 
	$$
	For $l+1 \leq v \leq 2l+1$ its degree $j$ component is given by
	\eqref{cocycle}
	while for $2l+2 \leq v \leq k$ it is given by 
	$$
	c(k,l)_v ((\varphi, e) \wedge (\widetilde{\varphi}, \widetilde{e})) = 
	\sum_{p = 0}^{l}  \Big(\varphi_{j-p}(Tr(\widetilde{e}_p)) + 
	\widetilde{\varphi}_{j-p}(Tr(e_p))  \Big)
	$$
\end{proposition} 
 \textit{Proof.} To show that the cocycle is relative observe that the 
 second sum in \eqref{cocycle} does not involve $e_0$ as $v - l \geq 0$
 the first sum does involve terms like $[s\widetilde{\varphi}, e_0]$ but those
 only depend on the covariant derivative of $e_0$ which vanishes by assumption.
 $\square$

\subsection{Torsors, associated to order $k$ extensions.}

We continue with the notation in the introduction
For $k \geq 1$, we would like to find an order $k$ extension of $E$.  
Given such extension $E^{(k)}$, by 
\eqref{lemma: filtered k isomorphism} on any affine open set $U \subset X$ there exist isomorphisms 
of filtered sheaves of algebras and modules 
$$
\mathcal{O}^{(\infty)} \simeq 
\widehat{Sym}_{\mathcal{O}_X} N^*,  \qquad  E^{(k)} \simeq 
Sym^{\leq k}_{\mathcal{O}_X} N^*   \otimes_{\mathcal{O}_X} E 
$$
where $\widehat{Sym}$ stands for the completed symmetric power, i.e. the infinite product of
$Sym^i$ for all $i \geq 0$. In particular, the respective isomorphisms hold 
for stalks 
at any closed point $x \in X$. Completing at the maximal ideal $\mathfrak{m}_x$ of $x$
and denoting $E_{x,0} = E_x/\mathfrak{m}_x E_x$
we conclude that, possibly shrinking $U$ to find a flat connection on $E|_U$
we can find isomorphisms
\begin{equation}
	\label{isomQ}
	\widehat{\mathcal{O}}_x^{(\infty)} \simeq 
	\widehat{Sym}_{\widehat{\mathcal{O}}_{X, x}} \widehat{N}_x^*,  \qquad  \widehat{E}_x^{(k)} \simeq 
	(Sym^{\leq k}_{\widehat{\mathcal{O}}_{X, x}} \widehat{N}_x^* ) \otimes_{\kappa} 
	E_{x, 0} 
\end{equation}
In addition, 
we can also find isomorphisms between completed objects
\begin{equation}
	\label{isomR}
	\widehat{\mathcal{O}}_{X, x} \simeq \mathcal{A}, \qquad
	\widehat{N}^*_x  \simeq \mathcal{S}, \qquad E_{x, 0} \simeq \mathcal{E}_0 
\end{equation} 
Combining  \eqref{isomQ} and \eqref{isomR}, we get filtered isomorphisms
\begin{equation}
	\label{isomP}
	\widehat{\mathcal{O}}^{(\infty)}_x \simeq 
	\mathcal{A}_S, \qquad \widehat{E}^{(k)}_x \simeq  \E^{\leq k}
\end{equation} 
Working at first informally, we define three infinite  dimensional torsors 
$$
Q(k) \to X, \qquad R \to X, \qquad P(k) \to X
$$
such that the fiber over $x \in X$ consists of all filtered isomorphisms as in 
\eqref{isomQ}, \eqref{isomR}, \eqref{isomP}, respectively. We assume that
isomorphisms of algebras are multiplicative and unital, and that isomorphisms
of modules are compatible with the action of algebras. We also assume that in 
all three cases the isomorphisms are compatible with filtrations at the maximal
ideals, and that for \eqref{isomQ}, \eqref{isomP} the filtration by the ideal 
sheaf of $X$ maps isomorphically to the degree filtration on the symmetric 
algebra.  It follows from the definitons that
$$
\qquad P(k) = Q(k) \times_X R,
$$
By definition, the three (infinite dimensional) fiber bundles are torsors over 
the automorphism groups of the objects on the right hand side of 
\eqref{isomQ}, \eqref{isomR}, \eqref{isomP}, respectively. In the case of $Q(k)$ we 
are dealing with the sheaf of algebras, which is a twisted form of 
functions with values in a prounipotent group over $k$. We denote the three 
groups by $U(k)$, $F$ and $G(k)$, respectively. 

\bigskip
\noindent 
First we describe the pro-agebraic group $G(k)$ as the group of filtered 
automorphisms of the pair $(\mathcal{A}_{\mathcal{S}}, \E^{\leq k})$. 
We assume that the automorphims preserve the filtration by powers of
the ideal $\mathcal{S} \cdot \mathcal{A}_{\mathcal{S}}$, and that
the induced automorphisms of associated graded objects are compatible
with the filtration by powers of the maximal ideal $\mathfrak{m} 
\subset \mathcal{A}$.  The group $F$ corresponding to \eqref{isomR},
is the product of the automorphisms of $(\mathcal{A}, \mathcal{S})$
which are compatible with the filtration induced by $\mathfrak{m}$, 
by the general linear group $GL(\E_0)$. Finally, the normal pro-unipotent subgroup $U(k) \subset G(k)$ is formed by the 
filtered automorphisms which induce identity on the associated 
graded quotient of the filtration by powers of $\mathcal{S} \cdot 
\mathcal{A}_{\mathcal{S}}$. We have a semidirect decomposition 
$G(k) = U(k) \rtimes F$. 

We emphasize that $Q_k$ is not a $U_k$-torsor since $F$
is not a normal subgroup.  In fact, the automorphisms of the two 
objects on the right hand side of \eqref{isomQ} form a Zariski 
sheaf of pro-unipontent groups. We will not use this approach,
just viewing $Q(k)$ as $P(k)/F$. 
If $P(k)$ is trivialized on a Zariski open
subset then $\pi_k: Q(k) \to X$
is a projective limit of morphisms with affine fibers. 
Therefore we can apply Lemma \ref{torsor extension}.

\bigskip
\noindent 
A more rigorous way to define it is to consider the sheaves $Jets(\mathcal{O}^{(\infty)})$, 
$Jets(E^{(k)})$ of jets of sections. By definition, $Jets(\mathcal{O}^{(\infty)})$ is the 
completion of $\mathcal{O}_X \otimes_\kappa \mathcal{O}^{(\infty)}$ at the ideal 
sheaf given by the kernel of
$$
\mathcal{O}_X \otimes_\kappa \mathcal{O}^{(\infty)} \to 
\mathcal{O}_X \otimes_\kappa \mathcal{O}_X \to \mathcal{O}_X
$$
where the first arrow is the reduction of the second tensor factor modulo 
the ideal sheaf of $X$ in $Y$, and the second arrow is the product map. The 
$\mathcal{O}_X$-modules structure is induced by the action on the left factor. 
The jet bundle 
$Jets(E^{(k)})$ is obtained similarly by completing the product $\mathcal{O}_X
\otimes_\kappa E^{(k)}$ along the action of the same ideal sheaf. 

Then we would like to construct a pro-algebraic scheme $P(k) \to X$ such that
giving 
a morphism $S \to P(k)$ is equivalent to giving a morphism $\rho: 
S \to X$ and filtered isomorphisms 
$$
\rho^* Jets(\mathcal{O}^{(\infty)}) \simeq \mathcal{O}_S \otimes_\kappa \mathcal{A}; 
\qquad 
\rho^* Jets(E^{(k)}) \simeq \mathcal{O}_S \otimes_\kappa \E^{\leq k}
$$
which are compatible with multiplicative structures. For trivial $E$
and the diagonal embedding $X \to Y = X \times X$ such a torsor in constructed
in Section 6 of \cite{VdB} and that construction may be repeated almost 
verbatim (the only extra ingredient is
Lemma \ref{lemma: filtered k isomorphism}) 
in our case as well. Same applies to the 
case of $Q(k)$ and $R$.

\bigskip
\noindent
\textbf{Remark.} In fact, $P(k)$ and $R$ have a richer structure of 
transitive
Harish-Chandra torsors in the sense explained in \cite{BK}. 
This means that the actions of Lie algebras of $G(k)$ and $F$, respectively, 
by the vector fields on the total space of the corresponding torsor, 
can be extended to a larger algebras of all derivations, 
not necessarily preserving the subspace defined by the maximal ideal 
at each point. Once standard consequence of having such a structure is that
every associated vector bundle of such a torsor automatically has a flat connection.
We will only use this when proving that $P(k)$ allows to recover 
$E^{(k)}$ by taking flat sections of the associated vector bundles 
that happens to be isomorphic to $Jets(E^{(k)})$. 

\subsection{Obstructions in the Lie cohomology model.}
If an order $l$ extension $E^{(l)}$ admits a further 
order $k$ extension $E^{(k)}$, we have  lifts of torsors
$$
P(k) \to P(l) \to X; \qquad Q(k) \to Q(l) \to X. 
$$
We want to stress that each individual component $c(l,k)_v$ is not
a cocycle as action of $Der_0'$ on $End^{l,k}_{ab}$ does not preserve
the grading by symmetric powers of $\mathcal{S}$. Consequently, 
the resulting obstruction classes on $X$ will not take values in 
a direct sum of $Sym^v \con \otimes End(E)$ but rather in a twisted 
version that we define now. For $l+1 \leq k \leq 2l+2$ define
\begin{equation}
\label{bundle-small}
\mathcal{E}nd^{k, l}:=  (J^{l+1}/ J^{k+1})
  \otimes_{\mathcal{O}^{(l)}}  End_{\mathcal{O}^{(l)}} (E^{(l)})
\end{equation}
where $J \subset \mathcal{O}^{(\infty)}$ is the completion of the ideal sheaf $I_X
\subset \mathcal{O}_X$. Note that $k-l \leq l+1$ imples that the tensor 
factor on the left is indeed a module over $\mathcal{O}^{(l)}$.
For $k \geq 2l+2$ such a definition would not work and in this case we
define 
\begin{equation}
\label{bundle-large}
\mathcal{E}nd^{k, l}_{ab}: = Ker(\mathcal{E}nd^{2l+1, l} \oplus \mathcal{O}^{(k)}
\to \mathcal{O}^{2l+1})
\end{equation}
where the arrow is the difference of the trace map and the canonical 
surjection.

\begin{proposition}
    For $k > l$, the pullback of $\mathcal{E}nd^{k, l}_{ab}$ to the total space of 
    $P(l)$ is isomorphic to the trivial bundle obtained from the $Der_k$-module 
    $End^{l, k}_{ab}$.
	If a choice of $E^{(l)}$ is fixed then existence
	of its degree $k$ extension 
	$E^{(k)}$ implies the vanishing of the class
	$$
	GF(P_l, c(l,k)) \in H^2(X, \mathcal{E}nd^{k, l}_{ab});
	$$
	with $c(l,k)$ as in Proposition \ref{proposition: cocycle}.
	If $l< k < 2l+2$ and 
	the vanishing occurs then the set of isomorphism classes of possible 
	extensions $E^{(k)}$ is nonempty and forms a torsor over 
	$	H^1(X, \mathcal{E}nd^{k, l})$.
\end{proposition}
\textit{Proof.} Indeed, after pullback to the total space of $P(l)$ we have
trivializations \eqref{isomP} which induce an isomorphism of the
pullback of $Jets(\mathcal{E}nd^{l, k}_{ab})$ and the $G(l)$-equivariant bundle
$\mathcal{O}_P \otimes_\kappa End^{l, k}_{ab}$. Since the 
de Rham cohomology of $Jets(\mathcal{E}nd^{l, k}_{ab})$ agrees with 
coherent cohomology of $\mathcal{E}nd^{l, k}_{ab})$, the statement 
about the class follows. 

If $E^{(k)}$ exists then the Gelfand-Fuks class in question vanishes by 
Lemma \ref{torsor extension}.

By the same Lemma \ref{torsor extension}, 
the vanishing of the obstruction class for $l+1 < k \leq 2l+1$ implies that
the torsor $P(l)$ may be lifted to a torsor $P(k)$. Now we finish as in 
Section 2 of \cite{BK}: the $(G(k), Der_k)$-module $\E^{\leq k}$ gives an associated 
bundle  on $X$ with a flat connection 
and $E^{(k)}$ can be recovered as a sheaf of its flat sections. 
$\square$

\bigskip
\noindent
\textbf{Examples.}

\begin{enumerate}
\item For $k = l+1$ the obstruction bundle $End^{l,k} \simeq Sym^{l+1} \con 
\otimes End(E)$ does not depend on the choice of $E^{(l)}$ although the 
cohomology class in general does. 

\item For $l = 0, k=1$ there are no commutators on the right hand side of 
\eqref{cocycle} and it can be re-written in terms of familiar classes. 
The Lie version of the Atiyah class is represented by the cocycle 
in $C^1(Der_l, Der_0'; \Omega^1_{\mathcal{A}/\kappa} 
\otimes_{\mathcal{A}} End_{\mathcal{A}}$ which sends $(\varphi, \psi)$
to $\nabla^\E(e_0)$. On the other hand, a derivation 
$\varphi_v: \mathcal{A} \to \mathcal{S}$ may be decomposed into 
a de Rham differential $\mathcal{A} \to \Omega^1_{\mathcal{A}/\kappa}$
and an $\mathcal{A}$-linear morphism $a_v: \Omega^1_{\mathcal{A}/\kappa} \to 
\mathcal{S}$. This gives a Lie algebraic version of the Kodaira-Spencer 
cocycle $KS \in C^1(Der_k, Der_0'; Hom_{\mathcal{A}}
(\Omega^1_{\mathcal{A}/\kappa}, \mathcal{S}))$. The Gelfand-Fuks classes of 
these cocycles are the usual Atiyah class $At(E) \in H^1(X, \Omega^1 \otimes
End(E))$ and the Kodaira-Spencer class $KS(X, Y) \in H^1(X, Hom(\Omega^1, \con))$.
The obstruction class is simply their cup product, taking values in 
$H^2(X, \con \otimes End(E)$. This first order obstruction formula is well-known, 
cf. \cite{HT}.

\item When $E$ has rank 1, $End^{k, l}_{ab} \simeq J^{l+1}/J^{k+1} \simeq I_X^{l+1}
/I_X^{k+1}$.

\item If $E$ has rank 1 and $H^j(X, \sym^k \con) =0$
for $i=1,2$ and $k \geq 1$, we can conclude that there exists a unique
extension of every finite order $k$ (another known result, due to Grothendieck, 
cf. Proposition 3.12 in \cite{SGA2}). In some cases we can only guarantee the vanishing for $k \geq l+1$ with 
fixed $l$ and in that case each choice of $E^{(l)}$ admits a unique extension of
any order $k> l$.

\item When $X$ is a proper open subvariety in a smooth projective surface, the
second cohomology vanishes since $X$ has cohomological dimension 1. Thus, 
extensions $E^{(k)}$ always exist but they may not be unique. 
\end{enumerate}

\section{The Dolbeault model} 

In this section we assume that the base field $\kappa$ is the field $\mathbb{C}$ of 
complex numbers and that $X \to Y$ is an embedding of Kahler manifolds (of course, 
it suffices to have a Kahler metric on some open neighborhood of $X$ in $Y$). 
For any holomorphic  bundle $V$ on $X$ denote by $A^\bullet_X(V)$ the Dolbeault 
complex $(\Lambda^\bullet (\Omega^{0,1}_X) \otimes V, \overline{\partial})$ of $V$. We can also
define the Dolbeault algebra $A^\bullet_X(\mathcal{O}^{(\infty)})$ 
and a module $A^\bullet_X(E^{(k)}$ over it using the completion 
procedure explained in Section 4 of \cite{Yu}.

This section is based on a general idea of \cite{Ka} and \cite{CCT}: that the sheaf
of algebras $\mathcal{O}^{(\infty)}$ may be described by replacing 
$\prod_{u \geq 0} \sym^u \con$ by a quasi-isomorphic sheaf of dg algebras, then
deforming the differential, then passing to cohomology. One can develop 
a Cech version of this procedure using Thom-Whitney normalization instead 
of Cech complex, as done in \cite{CCT}, but the globalized Dolbeault version is a more
explicit. In the case of the diagonal embedding $X \to X \times X$ it is 
due to Kapranov.

\subsection{Yu's Dolbeault model of the infinitesimal neighborhood.}

\begin{proposition}
\label{Yu model}
	The choice of the Kaehler metric induces an isomorphism of filtered algebras
	$$
	A^\bullet_X(\mathcal{O}^{(\infty)}) \simeq A^\bullet_X(\widehat{Sym} N^*)
	$$
	which does not agree with the differentials. The Dolbeault 
	differential $\overline{\partial}$ on the left hand
	side, when transferred to the right hand side via the above isomorphism, 
	takes the form $\overline{\partial} + \mathfrak{D}$ were 
	$\mathfrak{D} = \mathfrak{D}_1 + \mathfrak{D}_2 + \ldots$ is an 
	algebra derivation, with derivations 
	$$
	\mathfrak{D}_v: A^\bullet_X(\widehat{Sym}^\bullet \con) \to 
	A^{\bullet+1}_X(\widehat{Sym}^{\bullet + v} con)
	$$
	obtained by the Leibniz rule extension from first order operators 
	$$
	\mathfrak{A}_v : A^\bullet_X(\mathcal{O}_X) \to A^{\bullet +1}_X (Sym^v \con), 
	\qquad \mathfrak{L}_v: A^\bullet_X(\con) 
	\to A^{\bullet+1}_X(Sym^{v+1} \con), \qquad i \geq 1
	$$	
\end{proposition}
The explicit formulas for $\mathfrak{A}_v$, $\mathfrak{L}_v$ 
in terms of the metric connection of a 
Kahler metric in the neighborhood of $X \subset Y$, may be found in 
Theorem 5 of \cite{Yu}. In fact, a choice of Kahler metric may be replaced by 
a choice of weaker structure. 
We only note here that $\mathfrak{A}_i$ is a composition of 
a holomorphic de Rham differential $\partial$ and a product with a smooth section 
$\mathfrak{a}_v \in A^1_X(\Omega^1_X, Sym^v \con)$, and that the first order operator 
$\mathfrak{L}_v$ has symbol given by (the symmetrization of) $a_v$, as follows immediately 
from the Leibniz rule. 

\bigskip
\noindent
\textbf{Remarks.}

\begin{enumerate}

\item The operators $\mathfrak{A}_v$ and $\mathfrak{L}_v$ are adjoint 
to the anchor maps and higher brackets, respectively, of the $L_\infty$
algebroid structure on the homological shift $N[-1]$ (realized by shifting
its Dolbeault complex), as studied in \cite{Yu}. We will not use this language
although everything in our paper may be restated in such terms. 

\item In an earlier work, cf. \cite{Ka}, Kapranov considered an important special case of the
diagonal embedding $X \to X \times X$. In this case the situation simplifies
greatly: the operators $L_i$ are obtained by symmetrizing the higher covariant
derivatives of the curvature tensor. The "adjoint anchor maps" $A_i$ are all zero 
(this is obvious for $A_1$, for example). 

Note that in this case the question of extending a vector bundle $E$ on 
$X$ to 
the formal neighborhood of the diagonal is trivial: the pullback 
with respect to 
either projection $X \times X \to X$ is actually the
 extension to the whole of 
$X \times X$. In particular, all obstruction 
 classes vanish.  

\end{enumerate}

\subsection{Dolbeault model for extended bundles.}

Now assume we have fixed an order $k$ extension $E^{(k)}$ and
an isomorphism 
$A^\bullet_X(\mathcal{O}^{(\infty)}) \simeq A^\bullet_X(\widehat{Sym^\bullet} \con)$ as in Proposition \ref{Yu model}. 

\begin{proposition}
\label{our model}
There exists an isomorphism compatible with the filtered module 
structures 
$$
(A^\bullet_X(E^{(k)}), \overline{\partial}) 
 \simeq (A^\bullet_X(Sym^{\leq k} \con \otimes E), \overline{\partial} + 
 \mathfrak{D}^E) 
$$
where $\overline{\partial} + \mathfrak{D}^E$ is a module 
derivation compatible with the
algebra  derivation $\overline{\partial} + \mathfrak{D}$ from Proposition 
\ref{Yu model}
and $\mathfrak{D}^E = \mathfrak{D}_1^E + \mathfrak{D}_2^E + 
\ldots + \mathfrak{D}_k^E$ where $\mathfrak{D}_v^E:
A^\bullet_X(\bigoplus_{u = 0}^{k-v} Sym^u \con \otimes E) \to 
A^{\bullet + 1}_X(\bigoplus_{u = 0}^{k-v} Sym^{u+v} \con \otimes E) $
is obtained by Leibniz rule extension from a first order operator 	
$$
\mathfrak{M}_v: A^\bullet_X(E) \to A^{\bullet + 1}_X(Sym^v \con \otimes E) 
$$
with the scalar symbol $\mathfrak{a}_i$. If an order $k$ extension 
$E^{(k)}$ restricts to an order $l$ extension $E^{(l)}$ then 
any choice of such isomorphism for $E^{(l)}$ may be extended to $E^{(k)}$
without changing the operators $\mathfrak{M}_1, \ldots, \mathfrak{M}_l$. 
\end{proposition}
\textit{Proof.} As $\mathcal{O}^{(k)}$ is identified
(by a $C^\infty$ isomorphism) with $Sym^{\leq k}
\con$, this allows to view $E^{(k)}$ as a locally free sheaf of finite rank 
over $\mathcal{O}_X$. We can choose a $C^\infty$-splitting 
$$
E^{(k)} \simeq E \oplus I_X E^{(k)}
$$
(for example by taking an orthogonal complement to the second 
factor, with respect to any Hermitian metric). Thus 
$Sym^v \con \cdot E$ is well-defined as a sub-bundle of $E^{(k)}$ for 
$v = 1, \ldots, k$, and we obtain a required direct sum splitting. 
In particular, we can transfer the differential $\overline{\partial}$ 
from $A^\bullet_X(E^{(k)})$ to $A^\bullet_X(Sym^{\leq k} \con
\otimes E)$
and it will be automatically a module derivation compatible from the
differential transferred from $A^\bullet_X(\mathcal{O}^{(k)})$. 

By the derivation property, it is uniquely determined by its restriction 
to $E$ where it splits into components $\mathfrak{M}_v: A^\bullet_X(E) 
\to  A_X^{\bullet + 1}(Sym^v \con \otimes E)$ for $v = 0, \ldots, k$.
Since the isomorphism of bundles restricts to identity on the 
associated graded of the filtration, we see that $E_0 = \overline{\partial}$. 
The symbol property is an immediate consequence of the derivation 
property, and the extension assertion is a direct consequence of the 
recursive nature of the above proof. 
$\square$

\subsection{Construction of extensions and obstructions.} 
Conversely, assume that we have the operators $\mathfrak{M}_v, v = 1, 
\ldots k$ such 
that the corresponding operator $\overline{\partial} + \mathfrak{D}^E$
squares to zero and is multiplicatively compatible with 
$\overline{\partial} + \mathfrak{D}$. 
Then the kernel of this operator is a module over 
$Ker(\widehat{\partial} + \mathfrak{D}) = \mathcal{O}^{(\infty)}$ 
due to the multiplicative agreement of the two differentials. 
This kernel also have a filtration by  submodules with associated 
graded $Sym^{\leq k} \con \otimes E$. It follows that $E^{(k)} := Ker 
(\overline{\partial} + \mathfrak{D}^E)$  is projective over 
$\mathcal{O}^{(k)}$ and thus gives an order $k$ extension.

\bigskip
\noindent 
To construct the operators $E_i$ we fix an Hermitian metric 
$h$ on $E$ and consider the unique Chern connection 
$\nabla^E = \nabla + \overline{\partial}$ 
compatible with $h$, where $\nabla$ is the $(1,0)$ part of the
connection. Its
curvature $R_\nabla = \overline{\partial} \nabla + 
\nabla \overline{\partial}$
is a $(1, 1)$-form twisted by endomorphisms of $E$.
We now set 
$$
\mathfrak{M}_v = \mathfrak{a}_v \nabla + \mathfrak{m}_v 
$$
where $\mathfrak{m}_v \in A^1_X(Sym^v \con \otimes End(E))$. 

\bigskip 
\noindent 
We define $\mathfrak{Der}_k$  to be the algebra of derivations 
of the pair $(A^\bullet_X(\widehat{Sym}^\bullet \con), 
A^\bullet_X(Sym^{\leq k} \otimes \con))$ that have the
form $(\mathfrak{D}, \mathfrak{D}^E)$
as in Propositions \ref{Yu model} and 
\ref{our model}, except $\mathfrak{a}_v$ and $\mathfrak{m}_v$
take values in $A^\bullet_X(\ldots)$ rather than 
$A^1(\ldots)$. The differential on $\mathfrak{D}_k$ is 
given by $[\overline{\partial}, \ldots]$. 
The homological grading of $\mathfrak{D}er_k$ is inherited from that
of $A^\bullet_X$

As we require that $(\overline{\partial} + \mathfrak{D}^E)^2 = 0$, the 
problem of lifting $E^{(l)}$ to $E^{(k)}$ is equivalent to the problem 
of lifting a Maurer-Cartan element from $\mathfrak{D}_l$ to 
$\mathfrak{D}_k$. Hence we can apply the formalism of Section 2.

\bigskip
\noindent 
Now we assume that $k \leq 2l +1$. Then there is an abelian extension 
$$
0 \to A^\bullet_X(\bigoplus_{u = l+1}^k Sym^u \con \otimes  End(E))
\to \mathfrak{Der}_{k} \to \mathfrak{Der}_{l} \to 0.
$$
As in the case of Lie cohomology model, it admits a section where 
the missing operators $\mathfrak{M}_{l+1}, \ldots, \mathfrak{M}_{k}$
get filled in by the operators $(s\mathfrak{D}_{l+1}, \ldots, 
s\mathfrak{D}_{k})$, obtained by the Leibniz rule extension of 
$\mathfrak{a}_v\nabla $ to an operator 
$$
s\mathfrak{D}_v:
A^\bullet_X(\bigoplus_{u = 0}^{k-v} Sym^u \con \otimes E) \to 
A^{\bullet + 1}_X(\bigoplus_{u = 0}^{k-v} Sym^{u+v} \con \otimes E)
$$
This section is not compatible with the differential or the bracket.
This means that $s\varphi$ may not be a Maurer-Cartan
solution and we would like to correct this by adding nontrivial 
terms $m_{l+1}, \ldots, m_k$, respectively. As in Section 2, this 
gives an obstruction class. 

This obstruction class is a degree two element of 
$A^\bullet_X(\bigoplus_{u = l+1}^k Sym^u \con \otimes  End(E))$, but 
the Dolbeault differential is deformed by the differential $\mathfrak{D}^{End}$
induced by commutators with $\mathfrak{M}_v, v = 1, \ldots_l$. 
The following lemma is
established by direct computation
\begin{lemma}
For $l+1 \leq k \leq 2l+1$ the isomorphisms of 
 Propositions \ref{Yu model} and 
\ref{our model} induce an isomorphism 
$$
(A^\bullet_X(\mathcal{E}nd^{l,k}), \overline{\partial}) 
 \simeq (A^\bullet_X(\bigoplus_{u = l+1}^k Sym^u \con \otimes  End(E)), \overline{\partial} + 
 \mathfrak{D}^{End}) 
$$
A similar statement holds for $\mathcal{E}nd^{l,k}_{ab}$ if $k \geq 2l+1$.
\end{lemma}
As mentioned before, when we lift a Maurer-Cartan solution 
$\psi_l \in \mathfrak{D}er_l$ to $s\psi_l \in \mathfrak{D}er_k$, 
the Maurer-Cartan solution may fail to hold, as the lifting map 
$s$ does not agree with differential or bracket. 
Hence we would like to adjust $s \psi_l$ by adding
$\mathfrak{m}_v \in A^1_X(Sym^v \con \otimes E$ to $s\mathfrak{D}_v$
for $v = l+1, \ldots k$. Let $\psi_k$ be the resulting element 
\begin{proposition}
 The element $\psi_k$ satisfies the Maurer-Cartan solution if and 
 only if 
$$
- \overline{\partial} \mathfrak{m}_v - 
\sum_{p=l+1}^{v-1} [s\mathfrak{D}_{v-p}, 
\mathfrak{m}_p] =  \mathfrak{a}_v R_{\nabla} +
\frac{1}{2} \sum_{p=v-l}^l [\mathfrak{m}_{v-p}, 
\mathfrak{m}_p] +  \sum_{p=i-l}^l [s\mathfrak{D}_{i-p}, 
\mathfrak{m}_p]
$$
for $v = l+1, \ldots, k$, 
where $R_{\nabla} = \overline{\partial} \nabla + \nabla \overline{\partial}
\in A^1_X(\Omega^1_X \otimes  End(E))$
is the curvature of the connection $\nabla + \overline{\partial}$.
In other words, the right hand side should 
define a trivial cohomology class in $H^2(Y^{(l)}, \mathcal{E}nd^{l,k})$.
\end{proposition}
\textit{Proof.} We apply last equation in Section 2.2. The claimed formula 
folows once we prove that 
$$
\Delta_2(s\mathfrak{D}_p,s \mathfrak{D}_q) = 0, \qquad 
\Delta_1(s\mathfrak{D}_v) = - \mathfrak{a}_v R_\nabla
$$
The first assertion is derived from the identity $\nabla^2 =0$
that holds for our metric connection. This is established by a tedious local
computation in which we first reduce to the case when $\con$ is trivial, then to 
the case when $\mathfrak{L}_i = \mathfrak{a}_i d$ and finally to the case $\con =
\mathcal{O}_X$. The expression $\Delta_1(s\mathfrak{D}_v)$ is just
$[\overline{\partial}, s\mathfrak{D}_v]$ and the Leibniz rule allows
to reduce to the computation of the restriction to $A^\bullet_X(E)$
where the assertion follows from $[\overline{\partial}, \mathfrak{a}_v] = 0$
(which is a consequence of $(\overline{\partial} + \mathfrak{D})^2 = 0$). 
$\square$

\begin{proposition}
For $l+1 \leq k \leq 2l +1$ an order $l$ extension $E^{(l)}$ with a fixed
Dolbealut model as in Proposition \ref{our model} can be further 
extended to $E^{(k)}$ if and only if 
$$
 \mathfrak{a}_i R_{\nabla} +
\frac{1}{2} \sum_{p=i-l}^l [\mathfrak{m}_{i-p}, 
\mathfrak{m}_p] +  \sum_{p=i-l}^l [s\mathfrak{D}_{i-p}, 
\mathfrak{m}_p]$$
defines a zero class in $H^2(Y^{(l)},\mathcal{E}nd^{k, l})$. 
  If that is the case, the set of isomorphism classes of all extensions
  $E^{(k)}$ is a torsor over $H^1(Y^{(l)}, \mathcal{E}nd^{k, l})$
\end{proposition}
\textit{Proof.} By Proposition \ref{our model} the Dolbeault model
for $E^{(l)}$ can be extended to a Dolbeault model for $E^{(k)}$. 
By the previous proposition, existence of $\mathfrak{m}_{l+1}, 
\ldots, \mathfrak{m}_k$ is equivalent to the vanishing of the cohomology 
class, as claimed. A different Dolbeault model for the same $E^{(k)}$
will adjust these operators by an operator in the image of the differential, 
while in general a choce of $\mathfrak{m}_v$ may be adjusted by an element
in the kernel of the differential.
Hence isomorphism classes of order $k$ extensions are in bijective 
correspondence with the first cohomology group, as claimed. 
$\square$

\bigskip
\noindent
For $k \geq 2k+2$ we proceed as before, establishing a necessary condition.
In this case, the obstruction is degree 2 element in 
$$
A^2_X(\bigoplus_{u = l+1}^{2l+1} Sym^u \con \otimes  End(E))
\oplus A^2_X(\bigoplus_{v = 2l+2}^k Sym^v \con)
$$
and the complex involved is isomorphic to the Dolbeault complex of
$End^{l,k}_{ab}$. The same argument as for $k \leq 2k+1$ gives
\begin{proposition}
If for $k > 2l +1$ an order $l$ extension $E^{(l)}$   can be further 
extended to $E^{(k)}$ then
$$
\bigoplus_{v=l+1}^k
 \big(\mathfrak{a}_i R_{\nabla} +
\frac{1}{2} \sum_{p=v-l}^l [\mathfrak{m}_{i-p}, 
\mathfrak{m}_p] +  \sum_{p=i-l}^l [s\mathfrak{D}_{v-p}, 
\mathfrak{m}_p]\big)
\oplus  
\bigoplus_{v= 2l+2}^k
\big(\mathfrak{a}_v Tr(R_\nabla) + \sum_{p=1???}^l 
\mathfrak{D}_{v-p} Tr(\mathfrak{m}_p)\big)
$$
defines a zero class in $H^2(Y^{(l)},\mathcal{E}nd^{k, l}_{ab})$. 
\end{proposition}

\bigskip
\noindent
\textbf{Remarks.}

\begin{enumerate}
\item Note that $Tr(R_\nabla) \in A^1_X(\Omega^1_X)$ represents the first Chern 
class of $E$

\item We could as whether the vanishing of the class in the previous propositon 
ensures existence of some geometric object. It certainly gives an 
order $(2l+1)$ extension of $E$, since we can simply truncate the statement.
But if we look closer, we see that it also gives an order $k$ extension of the
determinant line bundle $\Lambda^e E$. This is perhaps the best we can expect
from looking at the usual - abelian - cohomology groups. 
\end{enumerate}

\section{The Cech model}

In this section we consider the Cech model of the situation. We choose
an open covering $\{U_i\}_{i \in I}$ 
such that $\mathcal{O}^{(\infty)}$ can be identified 
with the completed symmetric algebra of the conormal bundle on each 
$U_i$ (by Lemma \ref{lemma: filtered k isomorphism} 
this works for any affine covering) and 
on which $E$ admits a flat connection (e.g. because $E$ is trivial on each 
$U_i$, which holds for $U_i$ small enough). 

Unlike in the previous two models, we only consider extensions of first and
second order here. There reason is that, for a sheaf of Lie algebras $\mathcal{L}$
on $X$,  the Cech complex of the on the covering with coefficients in $\mathcal{L}$
does not in general have the structure of a dg Lie algebra, only a homotopy 
Lie algebra, also called $L_\infty$-structure. An alternative would be to 
replace the Cech complex by Thom-Sullivan normalization but that makes
the whole statement less transparent and computable.

For extensions of second order we can reduce to the case of a homotopy Lie
algebra where all higher $n$-ary operations vanish for $n \geq 3$. 
In fact, for first order extensions even the Lie bracket is equal to zero. 

Although the general case is not that difficult, we leave it to a motivated
reader, observing here that a general group cocycle condition below is
equivalent to a Maurer-Cartan condition if we use the Baker-Campbell-Hausdorff
formula and then observe that its degree $n$ term will match precisely the 
degree $n$ term of the Maurer-Cartan equation. This nontrivial result is 
implicit in \cite{Ge}, and is proved explicitly in \cite{FMM}. Note that both the 
$L_\infty$-structure and the BCH formula are defined up to certain ambiguity 
(in the case of BCH formula this is a consequence of the Jacobi identity), but
the formalism developed in \cite{Ge} implies that a choice of Dupont homotopy 
allows to fix both the $L_\infty$-structure and the CBH formula in a coherent 
way, which ensures agreement of the Maurer-Cartan equation and the cocycle condition. 

\subsection{Derivations on double intersections.} 
 Fix an affine open covering $X = \bigcup U_i$ and filtered isomorphisms
 \begin{equation}
 	\widetilde{\Phi}_i :  \hato^{(\infty)}|_{U_i}  \rightarrow  \widehat{\sym}^\bullet 
 	\con|_{U_i} 
 \end{equation}
which restrict to identity on associated graded quotients. These exist by 
Lemma \ref{lemma: filtered k isomorphism}.
On a double intersection
$U_{ij} = U_i \cap U_j$ the disagreement of two isomorphisms  
$\Phi_{ij} = \widetilde{\Phi}_j \widetilde{\Phi}_i^{-1}$ is a filtered algebra automorphism of $\widehat{\sym}_{\ox} \con|_{U_i} $ which restricts to identity on associated graded quotients. Therefore its
logarithm is well defined on our completed objects and we can write 
$$
\Phi_{ij} = exp(\varphi_{ij}) 
$$
\begin{lemma} 
\label{cech-functions}
	The map 
	$$
	\varphi_{ij}: \widehat{\sym}_{\ox} \con|_{U_{ij}}
	\to  \widehat{\sym}_{\ox} \con|_{U_{ij}}
	$$
	is a filtered $\kappa$-algebra derivation which vanishes on associated graded quotients. 
	In particular, it is uniquely determined by its restrictions 
	$$
	\mathbf{A}_{ij}: \ox|_{U_{ij}} \to \widehat{\sym}^{\geq 1}_{\ox} \con|_{U_{ij}}, \qquad 
	\mathbf{L}_{ij}: \con|_{U_{ij}} \to \widehat{\sym}^{\geq 2}_{\ox} \con|_{U_{ij}}.
	\square
	$$
\end{lemma} 

\bigskip
\noindent 
 We split the above two operators into homogeneous components: 
 $$
 \mathbf{A}^s_{ij}: \ox|_{U_{ij}} \to \sym^s_{\ox} \con|_{U_{ij}}, \qquad 
 \mathbf{L}^{s-1}_{ij}: \con|_{U_{ij}} \to \sym^s_{\ox} \con|_{U_{ij}}.
 $$
 	The following lemma is proved by a straightforward computation 
 \begin{lemma}
 	Each $\mathbf{A}^s_{ij}$ is a derivation and can be represented as a 
 	composition of $d: \ox|_{U_{ij}} \to \Omega^1_X|_{U_{ij}}$ and an 
 	$\ox$-linear operator $\mathbf{a}^s_{ij}:  \Omega^1_X|_{U_{ij}} \to \sym^s_{\ox} \con|_{U_{ij}}$.
 	Each $\mathbf{L}^{s-1}_{ij}$ is a first order differential operator with the symbol 
 	$1_{\con} \cdot \mathbf{a}_{ij}^{s-1}$.  
 \end{lemma} 

\bigskip
\noindent
Similar arguments apply to extensions $E^{(k)}$. As before,
we can choose isomorphisms 
\begin{equation}
	\widetilde{\Psi}_i :  E^{(k)}|_{U_i} \rightarrow  \sym^{\leq k} \con \otimes E |_{U_i}
\end{equation}
which restrict to identity on associated graded quotients, then form automorphisms 
$\Psi_{ij} = \widetilde{\Psi}_j \widetilde{\Psi}_i^{-1}$ of the term on the right, restricted to $U_{ij}$. 
As in the case of functions, we have 
$$
\Psi_{ij} = exp(\psi_{ij}) 
$$
Again, the following lemma is a direct consequence of defintions. 
\begin{lemma} The map 
\label{cech-modules}
	$$
	\psi_{ij}: \sym^{\leq k} \con \otimes E |_{U_{ij}}
	\to \sym^{\leq k} \con \otimes E |_{U_{ij}}
	$$
	agrees with the natural filtration on the truncated symmetric algebra and induces
	the zero map on the associated graded quotient. For sections $f$, resp. $e$, of
	$\sym^{\leq k} \con $,  resp. 
	$\sym^{\leq k} \con \otimes E $ on $U_{ij}$, one has the 
	following module derivation property:
	$$
	\psi_{ij} (f \cdot e) = \varphi_{ij}(f) \cdot e + 
	f \cdot \psi_{ij}(e) 
	$$
	In particular, $\psi_{ij}$ is uniquely determined by the restriction 
	$$
	\boldsymbol{M}_{ij}: E \to \bigoplus_{u = 1}^{k}\sym^{u} \con \otimes E
	$$
	Each component $\boldsymbol{M}^u_{ij}: E \to \sym^{u} \con \otimes E$
	is a first order operator with the symbol $1_E \otimes \mathbf{a}^u_{ij}$. For a connection 
	$\nabla_{i}$ on $E|_{U_i}$ one can write 
	$$
	\boldsymbol{M}^v_{ij} = \mathbf{a}^v_{ij} \nabla_{j} + \mathbf{m}^v_{ij}
	$$
	where $\mathbf{m}^v_{ij}: E \to \sym^{v} \con \otimes E$ is an $\ox$-linear 
	operator. $\square$
\end{lemma} 
 
As in the Dolbeault case, the operator $a^v_{ij} \nabla_j$ admits a Leibniz
rule extension to a derivation of $\sym^{\leq k} \con \otimes E$, increasing
the symmetric power degree by $v$. We
denote this operator by $s \varphi^v_{ij}$. 
 
 \subsection{The group cocycle condition and the Maurer-Cartan equation.} 
 
In addition to the multiplicative properties of our operators on each 
$U_{ij}$ we have the cocycle
condition on triple intersections $U_{ijh} = U_i \cap U_j \cap U_h$. We fix a total ordering on 
the index set of the covering and assume that $i < j < h$.  The group valued cocycle conditions
are, as usual, 
$$
exp(\varphi_{ih}) = exp(\varphi_{ij}) exp(\varphi_{jh}), \quad
exp(\psi_{ih}) = exp(\psi_{ij}) exp(\psi_{jh})
$$
where $(\varphi_{ij}, \psi_{ij})$ is considered a section over
$U_{ij}$ of the sheaf of Lie algebras 
$
\mathbf{Der}_k^+
$
formed by pair derivations $(\varphi, \psi)$ that satisfy the multiplicative 
conditions on Lemmas \ref{cech-functions} and 
\ref{cech-modules}, agree with filtrations and induce a zero map on associated 
graded quotients. 

If $\mathfrak{g}$ is the pro-nilpotent Lie algebra of a pro-unipotent algebraic group $G$, then 
for $x, y, z \in \mathfrak{g}$ then the group valued equation  $exp(z) = exp(x) exp(y)$ can be 
rewritten as 
$$
z = log(exp(x) exp(y)) =  x + y +  \frac{1}{2} [x, y] + \frac{1}{12}([x, [x, y]] - [y, [x, y]]) - \frac{1}{24} 
[y, [x, [x, y]]] + \ldots
$$
where the right hand side is given by the Baker-Campbell-Hausdorff formula, 
cf. \cite{Se}.

\bigskip
\noindent 
Note that for $\mathbf{Der}_k^+$ we also have the grading 
$\mathbf{Der}_k^+ = \prod_{t \geq 1} \mathbf{Der}_k^t$ where $\mathbf{Der}_k^t$ is formed by the 
elements $(\varphi^t, \psi^t)$ the increase the symmetric power degree 
in $\widehat{\sym} \con$ by $t$. 
As mentioned above, we are interested in the case when $t = 1, 2$.
Since the Lie brackets are compatible with this grading,
in degrees 1 and 2 the above cocycle conditions gives 
$$
\varphi^1_{ih} = \varphi^1_{ij} + \varphi^1_{jh} , \qquad 
\psi^1_{jh} = \psi_{ij}^1 + \psi^1_{jh} 
$$
$$
\varphi^2_{ih} = \varphi_{ij}^2 + \varphi^2_{jh} + \frac{1}{2}[\varphi^1_{ij}, \varphi^1_{jh}] , 
\qquad 
\psi^2_{ih} = \psi_{ij}^2 + \psi^2_{jh} + \frac{1}{2}[\psi^1_{ij}, \psi^1_{jh}] .
$$
We would like to rephrase these equations as cohomology
conditions on $\ox$-linear operators
$\mathbf{a}^1_{ij}$, $\mathbf{a}^2_{ij}$,  $\mathbf{m}^1_{ij}$, $\mathbf{m}^2_{ij}$.

\begin{proposition}
\label{cech obstruction}
A first order extension $E^{(1)}$ exists if an only if there exists
a 1-cochain $\mathbf{m}^1$ with values in $\con \otimes End(E)$ such that
$$
 - \delta \mathbf{m}^1 = \mathbf{a}^1_{ij} At(E)_{jh} 
$$
where $At(E)_{jh} = (\nabla_j - \nabla_h)$ is the Atiyah cocycle of $E$. 
Two choices of a resolving cochain $\mathbf{m}^1$
give isomorphic first order extensions if and only if their difference
is exact. 

\bigskip
\noindent
If a choice of $\mathbf{m}^1$ is fixed and the local connections $\nabla_1$
are flat, an extension of $E^{(1)}$ to $E^{(2)}$ exists if and only if there is 
a 1-cochain $\mathbf{m}^2$ with values in $\sym^2 \con \otimes E$ such that
$$
 - \delta \mathbf{m}^2 = \mathbf{a}^2_{ij} At(E)_{jh} + \frac{1}{2} [\mathbf{m}^1_{ij}, \mathbf{m}^1_{jk}]
 +  \frac{1}{2} [s\varphi^1_{jh}, \mathbf{m}^1_{ij}]
 + \frac{1}{2} [s\varphi^1_{ij}, \mathbf{m}^1_{jh}]
$$
In other words the right hand side defines the zero class in $H^2(X, 
\sym^2 \con \otimes E)$. Two choices of a resolving cochain $\mathbf{m}^2$
give isomorphic second order extensions if and only if their difference
is exact. 

\bigskip
\noindent
Finally, if any choice of $E^{(2)}$ exists (i.e. without fixing  $E^{(1)}$)
then $(\mathbf{a}^1_{ij} At(E)_{jh}, \mathbf{a}^2_{ij} Tr(At(E)_{jh}))$ defines
a zero cohomology class in $H^2(Y^{2}, \mathcal{E}nd^{0,2}_{ab})$
\end{proposition}
We recall that the obstruction bundle $\mathcal{E}nd^{0,2}_{ab}$ is
a non-trivial extension 
$$
0 \to \sym^2 \con \to \mathcal{E}nd^{0,2}_{ab} \to \con \otimes End(E) \to 0
$$
of sheaves of $\mathcal{O}^{(2)}$-modules (in fact, even $\mathcal{O}^{(1)}$-modules). We also note that $Tr(At(E)_{ij}$ represents the first
Chern class of $E$ and, while $\mathbf{a}^2$ is not a cocycle, the product
$\mathbf{a}^2_{ij} Tr(At(E)_{jh})$ can be seen to be a cocycle once we know that 
the class of $\mathbf{a}^1_{ij} At(E)_{jh}$ is trivial. 

\bigskip
\noindent 
\textit{Proof.}
By Leibniz rule the cocycle condition on $\psi^1$ reduces to the 
case of $\mathbf{M}^1_{ij}$ (once we know it for $\varphi^1$). Hence
$$
- (\mathbf{m}^1_{ij} + \mathbf{m}^1_{jh} - \mathbf{m}^1_{ih})
= \mathbf{a}^1_{ij}\nabla_j +  \mathbf{a}^1_{jh}\nabla_h - 
\mathbf{a}^1_{jh}\nabla_h =  \mathbf{a}^1_{ij}(\nabla_j - \nabla_h) 
= \mathbf{a}^1_{ij} At(E)_{jh}
$$
since $\delta \mathbf{a}^1 =0$. Observe that, as in the Dolbeault setting, 
$a^1_{ij}$ is a Cech cocycle of the Kodaira-Spencer class (the extension 
class of the conormal short extact sequence) and $At(E)_{jh}$ represents
the Atiyah class of the bundle $E$. This defines $\mathbf{m}^1$ up to 
a cohomology class in $H^1(X, \con \otimes End(E))$. 

For the second order extension, we use Section 2.2. Then 
$\Delta_1(s\varphi^2)$ gives $\mathbf{a}^2_{ij} At(E)_{jh}$.
For a flat connection the terms $\Delta_2(s\varphi^1_{ij})$ vanish
and the remaining terms 
give the second formula claimed in the proposition. 
Finally, for $l=0$ and $k=2$ we rearrange the terms of the second order
equation. 
$$
 - \delta \mathbf{m}^2 - 
  \frac{1}{2} [s\varphi^1_{jh}, \mathbf{m}^1_{ij}]
 - \frac{1}{2} [s\varphi^1_{ij}, \mathbf{m}^1_{jh}]
 = \mathbf{a}^2_{ij} At(E)_{jh} + \frac{1}{2} [\mathbf{m}^1_{ij}, \mathbf{m}^1_{jk}]
$$
then apply traces, keeping in mind that $Tr([s\varphi^1, \mathbf{m}^1] = 
\varphi^1 (Tr (\mathbf{m}^1))$ and that the trace of the commutator is zero. 
Finally, we observe that this procedure will turn the left hand side 
into a complete differential of the second order element in the Cech complex for $\mathcal{E}nd^{0,2}_{ab}$. $\square$

\bigskip
\noindent
Email: vbaranov@uci.edu, hongseoc@uci.edu. 

\end{document}